\documentclass[11pt]{article}
\usepackage{epsfig,amssymb}
\usepackage{palatino}
\usepackage{verbatim}

\def \e {\MM{e}}

\def\Frac#1#2{ {\mbox{\large$\frac{#1}{#2}$}}}

\def \R  {{\mathbb{R}}}

\def \N   {{\mathbb{N}}}

\def \SS {{\cal S}}

\def \OO {{\cal O}}

\def \wh {\widehat}

\makeatletter \@addtoreset{equation}{section} \makeatother

\newtheorem {theorem}    {Theorem} [section]

\newtheorem {proposition}[theorem] {Proposition}
\newtheorem {lemma}      [theorem] {Lemma}
\newtheorem {conjecture} [theorem] {Conjecture}

\newtheorem {remark}     [theorem] {Remark}

\newtheorem {standard form}   [theorem] {Standard Form}

\def\be {\begin{equation}}
\def\ee {\end{equation}}
\def\ba {\begin{eqnarray}}
\def\ea {\end{eqnarray}}
\def\baa {\begin{eqnarray*}}
\def\eaa {\end{eqnarray*}}

\def \sign {{\rm sgn\,}}

\newcommand{\supp}{{\rm supp\,}}

\def\proof{\noindent{\bf Proof.} }
\def\qed{\hfill $\Box$}
\newcommand {\lb} {\label}
\newcommand {\rf[1]} {(\ref{#1})}

\def \e {{\epsilon}}

\begin{document}
\title
{On the exact constant in Jackson-Stechkin inequality \\
for the uniform metric}

\author{S.\,Foucart, Yu.\,Kryakin, and A.\,Shadrin}

\date{}
\maketitle

\begin{abstract}
The classical Jackson-Stechkin inequality estimates
the value of the best uniform approximation of a periodic function $f$ by
trigonometric polynomials of degree $\le n-1$ in terms of its $r$-th modulus
of smoothness $\omega_r(f,\delta)$.  It reads
$$
    E_{n-1}(f) \le c_r\, \omega_r\Big(f,\frac{2\pi}{n}\Big)\,,
$$
where $c_r$ is {\it some} constant that depends only on $r$.
It was known that $c_r$ admits the estimate $c_r < r^{ar}$
and, basically, nothing else could be said about it.

The main result of this paper is in establishing that
$$
     (1-\Frac{1}{r+1})\, \gamma_r^* \le c_r < 5\, \gamma_r^*, \qquad
\gamma_r^* = \frac{1}{ {r \choose \lfloor \frac{r}{2} \rfloor} }
\asymp \frac{r^{1/2}}{2^r}\,,
$$
i.e., that the Stechkin constant $c_r$, far from increasing with $r$,
does in fact decay exponentially fast.
We also show that the same upper bound is valid for the constant $c_{r,p}$
in the Stechkin inequality for $L_p$-metrics with $p \in [1,\infty)$,
and for small $r$ we present upper estimates which are sufficiently
close to $1\cdot \gamma_r^*$.
\end{abstract}



\section{Introduction}

The classical Jackson-Stechkin inequality estimates
the value of the best uniform approximation of a periodic function $f$ by
trigonometric polynomials of degree $\le n-1$ in terms of its $r$-th modulus
of smoothness $\omega_r(f,\delta)$.  It reads
\be \lb{S}
    E_{n-1}(f) \le c_r\, \omega_r\Big(f,\frac{2\pi}{n}\Big)\,,
\ee
where $c_r$ is {\it some} constant which depends only on $r$
(see \cite{s51} or \cite[p.205]{dl}).

Besides the case $r=1$, hardly any attempts have been made
to find the best value of this constant $c_r$,
or even to determine its dependence on $r$.
Stechkin's original proof \cite{s51}
(as well as alternative ones)
allows to obtain the estimate $c_r < r^{ar}$,
and, basically, nothing else could be said about it.

The main result of this paper is in establishing that
\be \lb{c_r}
     (1-\Frac{1}{r+1})\, \gamma_r^* \le c_r < 5\, \gamma_r^*, \qquad
\gamma_r^* = \frac{1}{ {r \choose \lfloor \frac{r}{2} \rfloor} }
\asymp \frac{r^{1/2}}{2^r}\,,
\ee
i.e., that the Stechkin constant $c_r$, far from increasing with $r$,
does in fact decay exponentially fast.

We also show that the same upper bound is valid for the constant $c_{r,p}$
in the Stechkin inequality for $L_p$-metrics with $p \in [1,\infty)$,
and for small $r$ we present upper estimates which are sufficiently
close to $1\cdot \gamma_r^*$.

In retrospect, such a result could have been anticipated, since for
trigonometric approximation in $L_2$-metric, already in 1967, Chernykh
\cite{c67b} established that
\be \lb{C}
    E_{n-1}(f)_2
\le c_{r,2}\, \omega_r\Big(f,\frac{2\pi}{n}\Big)_2\,, \qquad
    c_{r,2} = \frac{1}{\sqrt{{2r \choose r}}} \asymp \frac{r^{1/4}}{2^r}\,,
\ee
proving also that such a $c_{r,2}$ is best possible
(for the argument $\delta = \frac{2\pi}{n}$ in $\omega_r(f,\delta)$).
However, this result was based on Fourier technique for $L_2$-approximation
and that does not work in other $L_p$-metrics.

Our method of proving \rf[c_r]
is based on deriving first the intermediate inequality
\be \lb{B}
   \|f\| \le c_{n,r}(\delta)\, \omega_r(f,\delta)\,, \qquad
   f \in T_{n-1}^\perp,
\ee
which is valid for the functions $f$ which are orthogonal
to the trigonometric polynomials of degree $\le n-1$.
This may be viewed as a difference analogue
of the classical Bohr-Favard inequality for differentiable functions
$$
   \|f\| \le \frac{F_r}{n^r}\,\|f^{(r)}\|\,, \qquad
   f \in T_{n-1}^\perp,
$$
and is of independent interest.

We make a pass from the Bohr-Favard-type inequality \rf[B]
to the Stechkin one \rf[S] by approximating $f$ with
the de la Vall\'ee Poussin sums $v_{m,n}(f)$ and using the fact that
$$
   f - v_{m,n}(f) \in T_m^\perp\,, \quad
   \omega_r(f\!-\!v_{m,n}(f),\delta)
               \le (1+\|v_{m,n}\|)\, \omega_r(f,\delta).
$$
With that we arrive at the inequality
$$
    E_{n-1}(f)
\le \|f - v_{m,n}(f)\|
    \le c_{m,n,r}(\delta) \omega_r(f,\delta)\,,
$$
where we finally minimize the resulting constant over $m$,
for given $r$, $n$ and $\delta$.

\section{Results} \lb{1}

For a continuous $2\pi$-periodic function $f$, we denote by $E_{n-1}(f)$
the value of best approximation of $f$ by trigonometric polynomials
of degree $ \le n-1$ in the uniform norm,
$$
    E_{n-1}(f) := \inf_{\tau \in T_{n-1}} \|f - \tau\|\,,
$$
and by $\omega_r(f,\delta)$ its $r$-th modulus of smoothness with the
step $\delta$,
$$
    \omega_r(f,\delta)
:= \sup_{0 < h \le \delta}\|\Delta_h^r(f,\cdot)\|, \qquad
  \Delta_h^r(f,x) = \sum_{i=0}^r (-1)^i {r\choose i} f(x+ih)\,,
$$
where $\Delta_h^r(f,x)$ is the forward difference of order $r$ of $f$
at $x$ with the step $h$.

We will study the best constant $K_{n,r}(\delta)$
in the Stechkin inequality
$$
    E_{n-1}(f) \le K_{n,r}(\delta)\,\omega_r(f,\delta)\,,
$$
i.e., the quantity
$$
     K_{n,r}(\delta)
  := \sup_{f\in C}\, \frac{E_{n-1}(f)}{\omega_r(f,\delta)}\,,
$$
which depends on the given parameters $n,r \in \N$ and $\delta \in [0,2\pi]$.

In such a setting (which goes back to Korneichuk and Chernykh)
we may safely consider $\delta = \frac{\alpha\pi}{n}$
with some $\alpha$ not necessarily $1$ or $2$.
The choice of particular $\delta$'s can be motivated by two reasons:

1)  "nice" look and/or tradition: $\delta = \frac{\pi}{n}$, or
$\delta = \frac{2\pi}{n}$, or (why not) $\delta = \frac{1}{n}$, and alike;

2) "nice" result:
$$
       \sup_{f\in C}\, \frac{E_{n-1}(f)}{\omega_r(f,\delta)}
\asymp c_{n,r}(\delta)\,.
$$
Ideally, both approaches should be combined to provide nice results
for nice $\delta$'s, but that happens not very often.

\medskip
In this paper we obtain the following results.

\medskip
1) First of all, we show that the exact order of the Stechkin constant
$K_{n,r}(\delta)$ at $\delta=\frac{2\pi}{n}$
(and in fact at any $\delta \in [\frac{2\pi}{n},\frac{\pi}{r}]$)
is $r^{1/2}2^{-r}$, namely
$$
     K_{n,r}(\Frac{2\pi}{n}) \asymp \gamma_r^* \asymp
     \frac{r^{1/2}}{2^r}\,,
$$
where
$$
    \gamma_r^*
 = \frac{1}{{r \choose \lfloor \frac{r}{2} \rfloor}}
 = \left\{\begin{array}{cl}
   \frac{1}{{2k \choose k}}, & r=2k;\\
   \frac{1}{{2k-1 \choose k-1}}, & r=2k-1.
   \end{array} \right.
$$
Moreover, we locate the exact value of this constant within
quite a narrow interval.

\medskip\noindent
{\bf Theorem 1.} {\it
We have
$$
         c'_r\big(\Frac{2\pi}{n}\big)\,\gamma_{r}^*
 \;\le\; \sup_{f\in C} \frac{E_{n-1}(f)}{\omega_{r}(f,\frac{2\pi}{n})}
 \;\le\; c_r\big(\Frac{2\pi}{n}\big)\,\gamma_{r}^*\,,
$$
where
$$
    c'_r\big(\Frac{2\pi}{n}\big)
  = \left\{ \begin{array}{cl}
    1 - \frac{1}{r+1}, & r = 2k-1; \\
    1,                 & r = 2k;
    \end{array} \right. \qquad n \ge 2r,
$$
and
$$
   c_r\big(\Frac{2\pi}{n}\big) = 5, \qquad n \ge 1.
$$
}
\quad
Surprising is the fact that, in this theorem, the upper estimate
is provided by one and the same linear method of approximation
that works for all $r$ simultaneously.
Namely, for any $r$, the de la Vall\'ee Poussin operator $v_{m,n}$
with $m = \lfloor\frac{8}{9}n\rfloor$ provides
$$
     \|f - v_{m,n}(f)\|
\le 5\,\gamma_r^* \omega_r\Big(f,\frac{2\pi}{n}\Big)\,, \qquad
     \forall r \in\N.
$$

2) Next, we show that the value of the constant $c_r(\delta)$
remains bounded uniformly in
$r$ and $n$ also for $\frac{\pi}{n} < \delta < \frac{2\pi}{n}$
(but it grows to infinity as $\delta$ approaches $\frac{\pi}{n}$).

\medskip\noindent
{\bf Theorem 2.} {\it
For any $\alpha > 1$, there
exists a constant $c_\alpha$ which depends only on $\alpha$ such that
$$
    E_{n-1}(f)
\le c_\alpha\,\gamma_{r}^*\,
    \omega_{r}\Big(f,\frac{\alpha\pi}{n}\Big)\,,\qquad
    n \ge 1\,.
$$
}

3) Thirdly, although we did not succeed
to reach the argument $\delta = \frac{\pi}{n}$ with an absolute constant
in front of $\gamma_r^*\, \omega_r(f,\delta)$, we prove that this constant
grows like $\OO(\sqrt{r}\ln{r})$ at most.

\medskip\noindent
{\bf Theorem 3.} {\it
For $\delta = \frac{\pi}{n}$, we have the estimate
$$
    E_{n-1}(f)
\le c_r(\Frac{\pi}{n})\,\gamma_{r}^*\,
    \omega_{r}\Big(f,\frac{\pi}{n}\Big)\,,\qquad
    c_r(\Frac{\pi}{n}) = \OO(\sqrt{r}\ln{r})\,, \qquad n \ge 1\,.
$$
}

4) Fourthly, for small $r$, the general upper bound
$c_r(\frac{2\pi}{n}) = 5$ can be decreased to the values that are
quite close to the lower bound $c_r' \approx 1$, thus giving support
to the (upcoming) conjecture that
$K_{n,r}(\delta) \le 1\cdot\gamma_r^*$ for $\delta \ge \frac{\pi}{n}$.

\medskip\noindent
{\bf Theorem 4.}  {\it
For $\delta = \frac{\pi}{n}$ and $\delta = \frac{2\pi}{n}$, we have
$$
    E_{n-1}(f)
\le c_{r}(\delta)\,\gamma_{r}^*\, \omega_{r}(f,\delta)\,,
$$
where $c_{2k-1}(\delta) = c_{2k}(\delta)$, and the values of
$c_{2k}(\delta)$ are given below
$$
\begin{array}{c|c}
c_2(\frac{\pi}{n}) & c_4(\frac{\pi}{n}) \\  \hline
1\frac{1}{4}       & 2\frac{7}{10}
\end{array}\,,
\qquad
\begin{array}{c|c|c}
c_2(\frac{2\pi}{n}) & c_4(\frac{2\pi}{n}) & c_6(\frac{2\pi}{n}) \\  \hline
1\frac{1}{16}       & 1\frac{1}{9}        & 1\frac{1}{2}
\end{array}\,.
$$
}

\medskip
5) Finally, all upper estimates in Theorems 1-4
remain valid for any $p \in [1,\infty]$.
(There is no need to give a separate proof of this statement,
since all the inequalities we used in the text still hold
for the $L_p$-metrics, $1 \le p < \infty$,
in particular the Bohr-Favard inequality \rf[BF] and the
inequalities of \S\ref{wv} involving the norms of
the de la Vall\'ee Poussin operator.)

\medskip\noindent
{\bf Theorem 5.} {\it
For any $p \in [1,\infty]$, we have
$$
    E_{n-1}(f)_p
\le c_{r}(\delta)\,\gamma_{r}^*\, \omega_{r}(f,\delta)_p\,,
$$
with the same constants $c_{r}(\delta)$ and the same $\delta$'s as
in Theorems 1--4. In particular,
$$
    E_{n-1}(f)_p
\le 5\,\gamma_{r}^*\, \omega_{r}\Big(f,\frac{2\pi}{n}\Big)_p\,.
$$
}

The latter $L_p$-estimate is hardly of the right order for
$1 < p < \infty$ because $\gamma_{r}^* \asymp r^{1/2}\,2^{-r}$,
while, for $p=2$, Chernykh's result \rf[C] says that
$K_{n,r}(\frac{2\pi}{n})_2 \asymp r^{1/4}\,2^{-r}$,
so one may guess that
$$
  K_{n,r}(\Frac{2\pi}{n})_p \asymp r^{\max(1/2p,1/2p')}\,2^{-r}\,.
$$
This guess is partially based on the results of Ivanov \cite{i94}
who obtained such an upper bound for the values $K_{n,r}(\delta)_p$
with relatively large $\delta = \frac{\pi r^{1/3}}{n}$, and proved that,
for $p \in [2,\infty]$, the order of the lower bounds is the same.

\medskip
6) The value $\delta = \frac{\pi}{n}$ is critical
in the sense that the Stechkin constant $K_{n,r}(\delta)$
and the constant $\gamma_r^*$ are no longer of the same (exponential)
order for $\delta = \frac{\alpha\pi}{n}$ with $\alpha < 1$.
Indeed, in this case, with $f_0(x) := \cos nx$,
we have $\omega_r(f_0,\delta) = 2^r \sin^r \frac{\alpha\pi}{2}$
(see \rf[cos]) and $E_{n-1}(f_0) = 1$,
so that, for $\alpha < 1$, we have
$$
    \frac{K_{n,r}(\Frac{\alpha\pi}{n})}{\gamma_r^*}
 >  \frac{cr^{1/2}}{\sin^r \frac{\alpha\pi}{2}}
 >  c_\alpha \lambda_\alpha^r, \qquad \lambda_\alpha > 1\,.
$$
This being said, a natural question arises from the two estimates
$$
     K_{n,r}(\Frac{2\pi}{n}) \asymp \gamma_r^*, \qquad
     K_{n,r}(\Frac{\pi}{n}) \le c \sqrt{r}\ln r\cdot \gamma_r^*
$$
whether an extra factor at $\delta = \frac{\pi}{n}$ is essential.
We believe it is not, and we are making the following brave conjecture.

\begin{conjecture} \lb{con}
For all $r\in\N$, we have
$$
    \sup_{n\in\N} K_{n,r}(\Frac{\pi}{n})
 := \sup_{n\in \N} \sup_{f\in C}\,
    \frac{E_{n-1}(f)}{\omega_r(f,\frac{\pi}{n})} = 1\cdot \gamma_r^*, \qquad
    \gamma_r^* = \frac{1}{{r \choose \lfloor \frac{r}{2} \rfloor}}\,.
$$
\end{conjecture}

(Our point is mainly about the upper bound, namely that
$K_{n,r}(\delta) \le 1\cdot\gamma_r^*$, for any $\delta \ge \frac{\pi}{n}$.
The lower bound for even $r=2k$ is established in this paper,
while for odd $r$ we guess that
$K_{n,r}(\delta)$ tends to  $\gamma_r^*$
at $\delta = \frac{\pi}{n}$ for large $n$,
but for $\delta > \frac{\pi}{n}$ it takes smaller values.)

This conjecture is true for $r=1$, for in this case we have
Korneichuk's result \cite{k62}:
$$
   1 - \frac{1}{2n} \le  K_{n,1}(\Frac{\pi}{n}) < 1.
$$

For $r=2$, the conjecture gives the estimate
$K_{n,2}(\frac{\pi}{n}) = \frac{1}{2}$
which is (to a certain extent) stronger than Korneichuk's one
(because $\omega_2(f,\delta) \le 2 \omega_1(f,\delta)$),
so it would be interesting to prove (or to disprove) it in
this particular case. Meanwhile, acccording to Theorems 1 and 4,
we have
$$
   \frac{1}{2}\le  K_{n,2}(\Frac{\pi}{n}) \le \frac{5}{8}\,,\qquad
   \frac{1}{2}\le  K_{n,2}(\Frac{2\pi}{n}) \le \frac{17}{32}\,.
$$

For arbitrary $r$, it seems unlikely that
the value of the Stechkin constant will ever be precisely determined,
but it would be a good achievement to narrow the interval for
$K_{n,r}(\frac{2\pi}{n})$, say, to $[\gamma_r^*,2\gamma_r^*]$,
and to settle down the
correct order of $K_{n,r}(\frac{\pi}{n})$ with respect to $r$.

7) We finish this section with the remark that if,
with some constant $c(\delta)$, the inequality
$$
   E_{n-1}(f) \le c(\delta)\, \gamma_{r}^*\, \omega_{r}(f,\delta)
$$
is true  for an even $r=2k$, then it is true for the odd $r=2k-1$ too,
with the same constant $c(\delta)$. Indeed, since
$\gamma_{2k-1}^* = 2 \gamma_{2k}^*$, and
$\omega_{2k}(f,\delta) \le 2\, \omega_{2k-1}(f,\delta)$, we have
\baa
    E_{n-1}(f)
&\le& c(\delta)\, \gamma_{2k}^*\, \omega_{2k}(f,\delta) \\
&\le& c(\delta)\, \gamma_{2k}^*\,\cdot 2\, \omega_{2k-1}(f,\delta)
 =    c(\delta)\, \gamma_{2k-1}^*\,\omega_{2k-1}(f,\delta)\,.
\eaa
Therefore, it is sufficient to prove upper estimates only for even
$r=2k$.


\section{Smoothing operators}

Here, we present the general idea of our method.

1) For a fixed $k$, with
$$
   \wh\Delta^{2k}_t(f,x)
:= \sum_{i=-k}^k (-1)^i {2k\choose k+i}f(x + it)
$$
being the central difference of order $2k$ with the step $t$, and with
$\phi_h$ being an integrable function which satisfies conditions
\be \lb{phi}
   a) \quad \phi_h(t) = \phi_h(-t), \qquad
   b) \quad \supp \phi_h = [-h,h], \qquad
   c) \quad \int_\R\phi_h(t)\,dt = 1,
\ee
consider the following operator
\be \lb{W_h}
      W_h(f,x)
 :=   \frac{1}{{2k \choose k}}
      \int_\R \wh\Delta^{2k}_t(f,x) \phi_h(t)\,dt\,.
\ee
If a given subspace $\SS$ is invariant under the operator $W_h$,
and if the restriction $W_h$ to $\SS$ has a bounded inverse,
then, for any $f\in\SS$, we have a trivial estimate
\be \lb{fW}
    \|f\| \le \|W_h^{-1}\|_\SS \|W_h(f)\|\,, \qquad f\in\SS.
\ee
It follows immediately from the definition that
\be \lb{Ww}
      \|W_h(f)\|
 \le \|\phi_h\|_1\,\gamma_{2k}^*\,\omega_{2k}(f,h), \qquad
    \gamma_{2k}^* = \frac{1}{{2k \choose k}}\,,
\ee
and we arrive at the following inequality:
\be \lb{cW}
    \|f\|
\le  c_{2k}(h)\,\gamma_{2k}^*\, \omega_{2k}(f,h)\,,
     \qquad c_{2k}(h) = \|\phi_h\|_1 \|W_h^{-1}\|_\SS\,,
\ee
valid for all functions $f$ from a given subspace $\SS$.

\medskip
2) Next, we present $W_h$ as $W_h = I-U_h$
what allows us to get some bounds for $\|W_h^{-1}\|$ in \rf[cW]
in terms of $U_h$.

To this end, for integer $i$ (and, in fact, for any $i$),
define the dilations $\phi_{ih}$ and the convolution operators $I_{ih}$
by the rule
$$
   \phi_{ih}(t) := \Frac{1}{i}\phi_h(\Frac{t}{i})\,, \qquad
   I_{ih}(f) := f*\phi_{ih} := \int_\R f(\cdot-t)\phi_{ih}(t)\,dt\,.
$$
Then, taking into account that
$$
   \int_\R f(x - it)\phi_h(t)\,dt
 = \int_\R f(x - \tau)\Frac{1}{i}\phi_h(\Frac{\tau}{i})\,d\tau
 = I_{ih}(f)\,,
$$
and that also $I_{ih} = I_{-ih}$ (because $\phi_{ih}$ is even),
we may put $W_h$ in the following form:
$$
    W_h
 = \frac{1}{{2k \choose k}}
       \sum_{i=-k}^{k} (-1)^{i} {2k \choose k+i} I_{ih} \\
 = I - 2 \sum_{i=1}^{k} (-1)^{i+1} a_i I_{ih}, \qquad
    a_i
 := \frac{{2k \choose k+i}}{{2k \choose k}}\,.
$$
So, with the further notations
$$
   U_{h} := 2 \sum_{i=1}^{k} (-1)^{i+1} a_i I_{ih}\,, \qquad
   \psi_{kh} := 2 \sum_{i=1}^{k} (-1)^{i+1} a_i \phi_{ih}\,,
$$
we obtain
$$
    W_h = I - U_{h}\,, \qquad U_{h}(f) = f*\psi_{kh}\,.
$$

Respectively, we may rewrite the inequality \rf[cW]
in the following way.

\begin{lemma}
If the opeartor $(I-U_h)^{-1}$ is bounded on a given subspace $\SS$,
then, for all $f\in\SS$, we have
$$
    \|f\|
\le  c_{2k}(h)\,\gamma_{2k}^*\, \omega_{2k}(f,h)\,,
     \qquad c_{2k}(h) = \|\phi_h\|_1 \|(I-U_{h})^{-1}\|_\SS\,.
$$
\end{lemma}

\medskip
3) Now, we call upon elementary properties of Banach algebras
(see, e.g., Kan\-to\-ro\-vich, Akilov \cite[Chapter 5, \S\,4]{ka})
for the claim that if an operator
$U: \SS \to \SS$ satisfies $\sum_{m=0}^\infty \|U^m\| < \infty$,
then the operator $I-U$ is invertible, and the norm of its inverse
admits the estimate
\be \lb{sum}
    \|(I-U)^{-1}\|_\SS \le \sum_{m=0}^\infty \|U^m\|_\SS\,.
\ee

\begin{proposition} \lb{rho}
If $\phi_h$ is such that
$\sum_{m=0}^\infty \|U_{h}^m\|_\SS = A_h < \infty$,
then, for any $f \in \SS$, we have
$$
    \|f\|
\le c_{2k}(h)\,\gamma_{2k}^*\, \omega_{2k}(f,h),
    \qquad  c_{2k}(h) = A_h\,\|\phi_h\|_1\,.
$$
\end{proposition}

\medskip
4) Finally, let us make a short remark about the structure of the
subspaces $\SS$ that may go into consideration.
It is clear that, if $\SS$ is shift-invariant, i.e., together
with $f$ it contains also $f(\cdot + t)$ for any $t$, then
$\SS$ is invariant under the action of $W_h$ for any $h$.
A typical example is a subspace $\SS$ that contains (or does not
contain) certain monomials $({\cos kx \atop \sin kx})$.

We will
consider $\SS = T_{n-1}^\perp$, the subspace of functions which
are orthogonal to trigonometric polynomials of degree $\le n-1$.


\section{A difference analogue of the Bohr-Favard inequality}

Denote by $T_{n-1}^\perp$ the set of functions $f$ which are orthogonal to
$T_{n-1}$, i.e., such that
$$
     \int_{-\pi}^\pi f(x) \tau(x)\,dx = 0, \qquad \forall \tau\in T_{n-1}\,.
$$
The Bohr-Favard inequality for such functions reads
\be \lb{BF}
   \|f\| \le \frac{F_r}{n^r} \|f^{(r)}\|, \qquad f \in T_{n-1}^\perp\,,
\ee
where $F_r$ are the Favard constants, which are usually defined by the formula
$$
   F_r
:= \frac{4}{\pi}\sum_{i=0}^\infty \frac{(-1)^{i(r+1)}}{(2i+1)^{r+1}}\,,
$$
and which satisfy the following relations:
$$
    F_0 = 1 < F_2 = \frac{\pi^2}{8} < \cdots < \frac{4}{\pi} < \cdots
    < F_3 = \frac{\pi^3}{24} < F_1 = \frac{\pi}{2}\,.
$$

In this section we obtain a difference analogue of the Bohr-Favard
inequality in the form
$$
   \|f\|
\le c_{n,2k}(h)\, \gamma_{2k}^*\,\omega_{2k}(f,h),
   \qquad f \in T_{n-1}^\perp\,,
$$
using the approach from the previous section (Proposition \ref{rho}).
Namely, we consider the operator
$$
    U_h = 2\sum_{i=1}^k (-1)^{i+1} a_i I_{ih}, \qquad
    \textstyle a_i = {2k \choose k+i}/{2k \choose k}\,,
$$
with the following specific choice of $I_h$ (and respectively of $\phi_h$):
$$
    I_h(f,x)
 := \frac{1}{h^2}\int_{-h/2}^{h/2}\int_{-h/2}^{h/2}f(x-t_1-t_2)\,dt_1\,dt_2\,,
$$
i.e., taking $I_h(f)$ as the Steklov function of order $2$.
It is known that  $I_{h}(f,x) = f*\phi_h$, where
$$
    \phi_h(t) = \left\{\begin{array}{ll}
    \frac{1}{h}(1-\frac{|t|}{h}), & t \in [-h,h], \\
      0, & \mbox{otherwise,}
    \end{array} \right.
$$
i.e., $\phi_{h}$ is the $L_1$-normalized B-spline of order $2$
(the hat-function) with the step-size $h$ supported on
$[-h,h]$. We also have
$$
   I_{ih}''(f,x)
 = - \frac{1}{(ih)^2}\wh\Delta^2_{ih}f(x)
 = \frac{1}{(ih)^2}[f(x-ih) - 2f(x) + f(x+ih)]\,.
$$
We denote by $\|U_h\|_{T_{n-1}^\perp}$ the norm of the operator
$U_h$ on the space $T_{n-1}^\perp$.

\begin{lemma} \lb{|U|}
We have
\be \lb{|U|e}
    \|U_h''\| \le \frac{\pi^2\mu^2}{h^2}\,, \qquad
    \mu^2 := \mu_{2k}^2
 := \frac{8}{\pi^2} \sum_{{\rm odd}\,i}^k \frac{a_i}{i^2} < 1\,.
\ee
\end{lemma}

\proof
1) We have
\ba
     U''_h(f,x)
&=&  2 \sum_{i=1}^{k} (-1)^{i+1}\, a_i I''_{ih}(f,x) \nonumber \\
&=&  2 \sum_{i=1}^{k} (-1)^{i+1}\, \frac{a_i}{(ih)^2}
      \Big[f(x-ih) - 2f(x) + f(x+ih)\Big] \nonumber \\
&=&  \frac{2}{h^2} \sum_{i=1}^{k} (-1)^{i+1}\, \frac{a_i}{i^2}
      \Big[f(x-ih) - 2f(x) + f(x+ih)\Big] \lb{f} \\
&=&  \frac{2}{h^2} \sum_{i=1}^{k} (-1)^{i+1} a_i'\, \Big[\!-\!2f(x)\Big]
     +  \frac{2}{h^2} \sum_{i=1}^{k} (-1)^{i+1}\,a_i'
      \Big[f(x-ih) + f(x+ih)\Big]\,, \nonumber
\ea
where in the last line we put $a_i' = \frac{a_i}{i^2}$. Hence,
\baa
       \frac{h^2}{4}\,\frac{\|U''_h(f)\|}{\|f\|}
&\le&  \Big|\sum_{i=1}^{k} (-1)^{i+1} a_i'\Big| + \sum_{i=1}^{k} |a_i'|
\;=\;  \sum_{i=1}^{k} (-1)^{i+1} a_i' + \sum_{i=1}^{k} a_i' \\
 &=&   2 \sum_{{\rm odd}\,i}^k a_i'
\;=\;   2 \sum_{{\rm odd}\,i}^k \frac{a_i}{i^2}
   =: \frac{\pi^2}{4}\,\mu^2
\eaa
i.e.,
$$
    \|U''_h\| \le \frac{\pi^2\mu^2}{h^2}\,.
$$

2) The estimate for $\mu^2$ follows from the fact that
$a_i = {2k \choose 2k+i}/{2k\choose k} < 1$, and
that $\sum_{i=1}^\infty \Frac{1}{i^2} = \Frac{\pi^2}{6}$:
$$
     \frac{\pi^2}{8} \mu^2
  =  \sum_{{\rm odd}\,i}^k \frac{a_i}{i^2}
  <  \sum_{{\rm odd}\,i}^\infty \frac{1}{i^2}
  =  \sum_{i=1}^\infty \frac{1}{i^2}
     - \sum_{{\rm even}\,i}^\infty \frac{1}{i^2}
  =  \Big(1-\frac{1}{4}\Big) \sum_{i=1}^\infty \frac{1}{i^2}
  =  \frac{\pi^2}{8}\,.
$$
We will prove in \S\ref{Spi/n} that
$1-\mu_{2k}^2 \asymp \frac{1}{\sqrt{2k}}$.
\qed

\begin{lemma} \lb{|U^m|}
We have
\be \lb{|U^m|e}
    \|U_h^m\|_{T_{n-1}^\perp}
\le F_{2m} \Big(\frac{\pi^2\mu^2}{n^2h^2}\Big)^m\,.
\ee
\end{lemma}

\proof
1) If $f$ is orthogonal to $T_{n-1}$, then so are its Steklov functions
$I_{ih}(f)$, hence $U_h(f)$ and the iterates $U_h^m(f)$ as well. Also,
the operators $D^2$ (of double differentiation) and $U_h$ commute
(since $D^2$ and $I_{ih}$ clearly commute).
Therefore, using the Bohr-Favard inequality with the $(2m)$-th derivative,
we obtain
\be \lb{BU}
    \|U_h^m(f) \|_{T_{n-1}^\perp}
\le \frac{F_{2m}}{n^{2m}}\, \|D^{2m} U_h^m(f)\|
 =  \frac{F_{2m}}{n^{2m}}\,\|\,[D^2 U_h]^m(f)\|
\le \frac{F_{2m}}{n^{2m}}\,\|D^2 U_h\|^m\|f\|\,.
\ee
By \rf[|U|e], we get $\|D^2 U_h\| \le \frac{\pi^2\mu^2}{h^2}$,
hence the conclusion.
\qed

\begin{remark} \rm
For $h = \frac{\pi}{n}$, we have equality in \rf[|U^m|e], i.e.,
$$
    \|U_h^m\|_{T_{n-1}^\perp} = F_{2m} \mu^{2m},
    \qquad h = \frac{\pi}{n}\,,
$$
which is attained on the Favard function $\varphi_n(x) = \sign\sin nx$.
Indeed, for $h = \frac{\pi}{n}$, we have
$$
    \varphi_n(x-ih) - 2 \varphi_n(x) + \varphi_n(x+ih)
= \left\{\begin{array}{cc}
  -4 \varphi_n(x), & {\rm odd}\;i, \\
   0, & {\rm even}\; i,
   \end{array} \right.
$$
and it follows from \rf[f] that
$U_h''(\varphi_n) = - \frac{\pi^2\mu^2}{h^2} \varphi_n$,
and respectively
$$
     D^{2m} U_h^m(\varphi_n)
  = (-1)^m \Big(\frac{\pi^2\mu^2}{h^2}\Big)^m \varphi_n\,.
$$
On the other hand, the Bohr-Favard inequality turns into equality
on the functions $f \in T_{n-1}^\perp$ such that
$f^{(2m)}(x) = a \varphi_n(x - b)$, hence on $U_h^m(\varphi_n)$.
Therefore, in \rf[BU], we have equalities all the way through.
\end{remark}

\begin{proposition} \lb{mu}
Let $f \in T_{n-1}^\perp$, and let $h > \frac{\pi}{n}\mu$. Then
\be \lb{gBh}
    \|f\|
\le c_{n,2k}(h)\,\gamma_{2k}^*\, \omega_{2k}(f,h)\,,
\ee
where
\be
   c_{n,2k}(h)
 = \Big(\cos \frac{\pi}{2}\rho\Big)^{-1}, \qquad
   \rho = \frac{\pi\mu}{nh} < 1.
\ee
\end{proposition}

\proof
From Proposition \ref{rho}, using the estimate \rf[|U^m|e], we obtain
$$
    c_{n,2k}(h)
 =  \sum_{m=0}^\infty \|U_h^m\|_{T_{n-1}^\perp}
\le \sum_{m=0}^\infty F_{2m}\rho^{2m}
 =  \Big(\cos \frac{\pi}{2}\rho\Big)^{-1}\,,
$$
the last equality (provided $\rho < 1$) being the
Taylor expansion of $\sec \frac{\pi}{2}x = 1/\cos \frac{\pi}{2}x$.
(The latter is usually given in terms of the Euler numbers $E_{2m}$ as
$\sec x = \sum_{m=0}^\infty \frac{|E_{2m}|}{(2m)!} x^{2m}$,
see, e.g., Gradshteyn, Ryzhik \cite[\S\,1.\,411.\,9]{gr}, so we have
$\sec \frac{\pi}{2}x
= \sum_{m=0}^\infty \frac{|E_{2m}|\pi^{2m}}{2^{2m}(2m)!} x^{2m}$,
and we use the fact that $F_{2m} = \frac{|E_{2m}|\pi^{2m}}{2^{2m}(2m)!}$,
see \cite[\S\,0.\,233.\,6]{gr}.)
\qed

\begin{theorem} \lb{gB}
If $f \in T_{n-1}^\perp$, then, for any $\alpha > 1$, we have
\be \lb{gBe}
    \|f\|
\le c_\alpha\,\gamma_{2k}^*\, \omega_{2k}(f,\Frac{\alpha \pi}{n})\,,
    \qquad
    c_\alpha
 =  \Big(\cos \frac{\pi}{2\alpha}\Big)^{-1}\,.
\ee
\end{theorem}

\proof
Just put $h = \frac{\alpha\pi}{n}$ in \rf[gBh], and use the fact that
$\mu < 1$.
\qed

\medskip
Let us give some particular cases of Theorem \ref{gB}.
\be \lb{aa}
\begin{array}{@{}ll@{}}
    1) \quad \alpha = 2,\qquad
    c_\alpha  = (\cos \frac{\pi}{4})^{-1} = \sqrt{2}, \qquad
&    \|f\| \le 1\frac{1}{2}\, \gamma_{2k}^*\,\omega_{2k}(f,\Frac{2\pi}{n})\,;
    \\[1ex]
    2) \quad \alpha = \Frac{3}{2},\qquad
    c_\alpha = (\cos \frac{\pi}{3})^{-1} = 2, \qquad
&    \|f\| \le 2\, \gamma_{2k}^*\,\omega_{2k}(f,\Frac{3\pi}{2n})\,;
    \\[1ex]
    3) \quad \alpha = \frac{4}{3},\qquad
    c_\alpha = (\cos \frac{3\pi}{8})^{-1} = 2.61, \qquad
&    \|f\| \le 2\frac{2}{3}\,\gamma_{2k}^*\,\omega_{2k}(f,\Frac{4\pi}{3n})\,;
    \\[1ex]
    4) \quad \alpha = \frac{5}{4},\qquad
    c_\alpha = (\cos \frac{2\pi}{5})^{-1} = 3.23, \qquad
&    \|f\| \le 3\frac{1}{4}\,\gamma_{2k}^*\,\omega_{2k}(f,\Frac{5\pi}{4n})\,.
\end{array}
\ee
From the relations
$
    \cos \frac{\pi}{2}x
  = \sin \frac{\pi}{2}(1-x) \ge \frac{\pi}{4}(1 - x^2),
$
it follows that, in \rf[gBe],
$$
    c_\alpha < \frac{4}{\pi} \Big(1 - \frac{1}{\alpha^2}\Big)^{-1}\,,
$$
i.e., $c_\alpha$ behaves like
$\frac{2}{\pi}\frac{1}{\alpha-1}$ as $\alpha \searrow 1$.

\begin{theorem}
If $f \in T_{n-1}^\perp$, then, for $\delta = \frac{\pi}{n}$, we have
\be \lb{gBmu}
    \|f\|
\le c_{2k}\,\gamma_{2k}^*\, \omega_{2k}(f,\Frac{\pi}{n})\,,
    \qquad
    c_{2k} = \OO(\sqrt{2k})\,.
\ee
\end{theorem}

\proof
Putting $h = \frac{\pi}{n}$ into \rf[gBh], we obtain the inequality
\rf[gBmu] with the constant
\be \lb{c_mu}
  c_{2k} =  \Big(\cos \frac{\pi}{2}\mu_{2k}\Big)^{-1}
 < \frac{4}{\pi}\Big(1 - \mu_{2k}^2\Big)^{-1}\,,
\ee
and we are proving in \S\ref{Spi/n} that
$1- \mu_{2k}^2 \asymp \frac{1}{\sqrt{2k}}$\,.
\qed


\section{Stechkin inequality for $\frac{\pi}{n} < \delta \le \frac{2\pi}{n}$}
\lb{wv}

1) Consider the de la Vall\'ee Poussin sum (operator)
\be \lb{v}
   v_{m,n} = \frac{1}{n-m}\sum_{i=m}^{n-1} s_i\,,
\ee
which is an average of $(n-m)$ Fourier sums $s_i$ of degree $i$.
For $m=n-1$ and for $m=0$, it becomes the Fourier sum $s_{n-1}$ and
the Fejer sum $\sigma_n = \frac{1}{n}\sum_{i=0}^{n-1} s_i$,
respectively.

Since $v_{m,n}(f)$ is the convolution of $f$ with the de la Vall\'ee Poussin
kernel $V_{m,n}$, we clearly have
$$
   \omega_k(v_{m,n}(f),\delta) \le \|v_{m,n}\|\, \omega_k(f,\delta),
$$
where $\|v_{m,n}\|$ is the norm, or the Lebesgue constant,
of the operator $v_{m,n}$.

Stechkin \cite{s51b} made a detailed studies of behaviour of the value
$\|v_{m,n}\|$ as a function of $m$ and $n$. We will need just two facts
from his work, one of them combined with a later result of Galkin
\cite{g71}.

a) The norm $\|v_{m,n}\|$ depends only on ratio $m/n$, and in a monotone
way. Precisely, with
$$
   \ell(x)
 := \frac{2}{\pi}\int_0^\infty \frac{|\sin xt\!\cdot\!\sin t|}{t^2}\,dt,
$$
which is (non-trivially) a monotonely increasing function of $x$, we have
$$
    \|v_{m,n}\| = \ell(x_{m/n}), \qquad x_{m/n} := \frac{1+m/n}{1-m/n}\,.
$$

b) The values of $\ell$ at integer points can be related to the
so-called Watson constants $L_{M/2}$ (for $M=2N$, they turn into the
Lebesgue constants $L_N := \|s_N\|$ of the Fourier operator $s_N$).
Namely,
$$
    \ell(M+1) = L_{M/2}\,,
$$
and from the result of Galkin \cite{g71} that
$L_{M/2} < \frac{4}{\pi^2}\ln(M+1) + 1$, we conclude that
\be \lb{p}
    \ell(p) < \frac{4}{\pi^2}\ln p + 1\qquad \mbox{for integer $p$},
\ee
therefore (rather roughly)
\be \lb{x}
    \ell(x) < \frac{4}{\pi^2}\ln (x+1) + 1 \qquad \mbox{for all $x$}.
\ee

\medskip
2) Now, from definition \rf[v], we see firstly that
$v_{m,n}(f)$ is a trigonometric polynomial of degree $\le n-1$, hence
$$
   E_{n-1}(f) \le \|f-v_{m,n}(f)\|,
$$
and secondly that $v_{m,n}$ acts as identity on $T_m$, therefore
$$
   f - v_{m,n}(f) \perp T_m\,.
$$
So, we may apply Proposition \ref{mu} to the difference  $f - v_{m,n}(f)$
to obtain
\baa
       E_{n-1}(f)
&\le& \|f - v_{m,n}(f)\| \\
&\le&  c_{m+1,2k}(h)\, \gamma_{2k}^*\,
       \omega_{2k}\Big(f-v_{m,n}(f),h\Big) \\
&\le& c_{m+1,2k}(h)\,(1 + \|v_{m,n}\|)\,\gamma_{2k}^*\,
      \omega_{2k}(f,h)\,\nonumber \\
& = & \Big[\cos \Big(\frac{\pi}{2}\frac{\pi\mu}{(m+1)h}\Big)\Big]^{-1}\,
      \Big[1 + \ell\Big(\frac{1+m/n}{1-m/n}\Big)\Big]\,
      \gamma_{2k}^*\, \omega_{2k}(f,h)\,.
\eaa
Now, with some parameter $s \in [0,1)$ which may well depend on
$n$ and $h$, we put in the last line
$$
   m = \lfloor sn \rfloor.
$$
With such an $m$, we have $m+1 > sn$ and $m/n \le s$, therefore
\be \lb{s1}
       E_{n-1}(f)
 \le  \Big[
      \cos \Big(\frac{\pi}{2} \frac{\mu}{s}\frac{\pi}{n h}\Big)\Big]^{-1}\,
      \Big[1 + \ell\Big(\frac{1+s}{1-s}\Big)\Big]
      \gamma_{2k}^*\,\omega_{2k}(f,h)\,.
\ee
Finally, taking $h = \frac{\alpha\pi}{n}$, and evaluating
the factor $1 + \ell(x_s)$ with the help of \rf[x], we obtain
\be \lb{s2}
       E_{n-1}(f)
 \le  \Big(\cos \frac{\pi\mu}{2\alpha s}\Big)^{-1}\,
      \Big[2 + \frac{4}{\pi^2}\ln\Big(\frac{2}{1-s}\Big)\Big]
      \gamma_{2k}^*\,\omega_{2k}\Big(f,\frac{\alpha\pi}{n}\Big)\,,
\ee
where we can minimize the right-hand side with respect to
$s \in (\frac{\mu}{\alpha},1)$.

\medskip
3) Now, using the last estimate, we establish Stechkin inequalities
for particular $\alpha$'s.

\begin{theorem} For all $n \ge 1$, we have
$$
       E_{n-1}(f)
  \le c\,\gamma_{2k}^*\,\omega_{2k}\Big(f,\frac{2\pi}{n}\Big)\,,
       \qquad c = 5\,.
$$
\end{theorem}

\proof
In \rf[s2], take $\alpha = 2$ and majorize $\mu$ by $1$.
Then the constant for $\delta=\frac{2\pi}{n}$ takes the form
$$
    c
=  \Big(\cos \frac{\pi}{4s}\Big)^{-1}\,
      \Big[2 + \frac{4}{\pi^2}\ln\Big(\frac{2}{1-s}\Big)\Big]\,.
$$
It turns out that the value $s = 8/9$ is almost optimal, and we obtain
Stechkin inequality with the constant
\be \lb{5}
   c = \Big(\cos \frac{9\pi}{32}\Big)^{-1}\,
      \Big[2 + \frac{4}{\pi^2}\ln 18\Big] = 4.999144 < 5.
\ee
To make sure that our step away from $5$ is free from a round-off error,
we notice that, for $s=\frac{8}{9}$, we have in \rf[s1]
$$
    \ell\Big(\frac{1+s}{1-s}\Big) = \ell(17) = L_8\,.
$$
Therefore, in the pass from \rf[s1] to \rf[s2],
we can use the estimate \rf[p] instead of \rf[x],
thus changing in \rf[5] the value $\ln 18$ to $\ln 17$,
and that will give the constant $c=4.962628$. We can make another bit down
by computing directly the Lebesgue constant $L_8 = 2.137730$, hence getting
$$
   c = \Big(\cos \frac{9\pi}{32}\Big)^{-1}\,
      \Big[1 + L_8\Big] = 4.946034,
$$
so that $c < 5$ is secured.
\qed

\begin{remark} \rm
Surprising is the fact that, in this theorem, the upper estimate
is provided by one and the same linear method of approximation
that works for all $r$ simultaneously.
Namely, for any $r$, the de la Vall\'ee Poussin operator $v_{m,n}$
with $m = \lfloor\frac{8}{9}n\rfloor$ provides
$$
     \|f - v_{m,n}(f)\|
\le 5\,\gamma_r^* \omega_r\Big(f,\frac{2\pi}{n}\Big)\,, \qquad
     \forall r \in\N.
$$
Perhaps it makes sense to try to derive such an estimate directly
from the properties of $v_{m,n}$.
\end{remark}

\begin{theorem} \lb{sa}
For any $\alpha > 1$, there exists a
constant $c_\alpha$ that depends only on $\alpha$ such that
\be \lb{al}
    E_{n-1}(f)
\le c_\alpha\,\gamma_{2k}^*\,
    \omega_{2k}\Big(f,\frac{\alpha\pi}{n}\Big), \qquad n \ge 1.
\ee
\end{theorem}

\proof
Putting (a non-optimal) $s = \frac{1}{\sqrt{\alpha}}$ in \rf[s2],
and again majorizing $\mu$ by $1$, we obtain \rf[al] with
\baa
      c_\alpha
 &=&  \Big(\cos \frac{\pi}{2\sqrt{\alpha}}\Big)^{-1}
      \left(\frac{4}{\pi^2}
      \ln\Big(\frac{2\sqrt{\alpha}}{\sqrt{\alpha}-1}\Big) + 2\right) \\
&\le& \frac{4}{\pi}\frac{\alpha}{\alpha-1}
      \left(\frac{4}{\pi^2}
      \ln\Big(\frac{2\sqrt{\alpha}}{\sqrt{\alpha}-1}\Big) + 2\right),
\eaa
where we have used the inequality $\cos\frac{\pi}{2}x \ge \frac{\pi}{4}(1-x^2)$
for $|x| \le 1$.
\qed


\section{Stechkin inequality for $\delta = \frac{\pi}{n}$} \lb{Spi/n}

\begin{theorem}
For $\delta = \frac{\pi}{n}$, and $r=2k$, we have
\be \lb{pi/n}
    E_{n-1}(f)
\le c_{r}(\Frac{\pi}{n})\,\gamma_{r}^*\,
    \omega_{r}\Big(f,\frac{\pi}{n}\Big)\,,\qquad
    n \ge 1,
\ee
where
\be \lb{cmu}
    c_{r}(\Frac{\pi}{n}) = \OO(\sqrt{r}\ln{r})\,.
\ee
\end{theorem}

\proof
From the estimate \rf[s2], with $h = \frac{\pi}{n}$ and $s=\sqrt{\mu}$,
we obtain the inequality \rf[pi/n] with the constant
\baa
     c_{2k}(\Frac{\pi}{n})
&=&  \Big(\cos \frac{\pi}{2}\sqrt{\mu}\Big)^{-1}
      \left(\frac{4}{\pi^2}
      \ln\Big(\frac{2}{1-\sqrt{\mu}}\Big) + 2\right) \\
&<&  \frac{4}{\pi}\frac{1}{1-\mu}\,
     \left(\frac{4}{\pi^2}
     \ln\Big(\frac{2}{1-\sqrt{\mu}}\Big) + 2\right).
\eaa
The estimate \rf[cmu] follows now from the fact that
$$1 - \mu_{2k}^2 > \frac{c_1}{\sqrt{2k}}\,, \qquad c_1 = \frac{2}{3},
$$
which we are proving in the next lemma.
With the value $c_1 = \frac{2}{3}$ at hands, we can give the explicit
estimate $c_{r}(\frac{\pi}{n}) < 2 \sqrt{r}\ln r + 12 \sqrt{r}$.
\qed

\begin{lemma}
For $\mu_{2k}^2
:= \frac{8}{\pi^2}\sum\limits_{{\rm odd}\,i}^k \frac{a_i}{i^2}$,
 where $a_i := {2k \choose k+i}/{2k\choose k}$, we have
\be \lb{mu<}
    \frac{c_1}{\sqrt{2k}} <  1-\mu_{2k}^2 <  \frac{c_2}{\sqrt{2k}}\,,
    \qquad c_1 = \frac{2}{3}, \quad c_2 = \frac{5}{4}\,.
\ee
\end{lemma}

\proof
Let us compute the value $\wh\Delta^{2k}_t(f_0,x)$ for $f_0(x) = \cos x$
at $x=0$. Since
$$
    \wh\Delta_t^2(\cos,x)
  = - \cos(x-t) + 2 \cos x - \cos(x+t)
  = 2 \cos x (1-\cos t) = 4 \sin^{2}\Frac{t}{2}\cos x\,,
$$
we have
$$
   \wh\Delta^{2k}_t(f_0,x)\Big|_{x=0}
 = 4^k \sin^{2k}\Frac{t}{2}\,.
$$
On the other hand, by the definition,
$$
   \wh\Delta^{2k}_t(f_0,x)\Big|_{x=0}
= \sum_{i=-k}^k (-1)^i {2k \choose k+i} \cos (x+it)\Big|_{x=0}
= {2k \choose k} \Big[1 - 2 \sum_{i=1}^k (-1)^{i+1} a_i \cos it\Big].
$$
So, we have
$$
    1 - 2 \sum_{i=1}^k (-1)^{i+1} a_i \cos it
  = \lambda_k \sin^{2k}\Frac{t}{2}, \qquad
    \lambda_k := \frac{4^k}{{2k\choose k}}\,.
$$
Integrating both parts twice, first time between $0$ and $u$, and then
between $0$ and $\pi$, we obtain: for the left-hand side
$$
    \Big[\frac{u^2}{2}
    + 2 \sum_{i=1}^k (-1)^{i+1} \frac{a_i}{i^2}\, \cos iu \Big]_0^\pi
 =  \frac{\pi^2}{2}
    - 4 \sum_{{\rm odd}\,i}^k \frac{a_i}{i^2}
 = \frac{\pi^2}{2}(1-\mu_{2k}^2)\,,
$$
and for the right-hand side
$$
    \lambda_k \int_0^\pi\int_0^u \sin^{2k}(\Frac{t}{2}) \,dt\,du
  = \lambda_k \int_0^\pi (\pi - t) \sin^{2k}(\Frac{t}{2}) \,dt
  = 4 \lambda_k \int_0^{\pi/2} \tau \cos^{2k}(\tau) \,d\tau
$$
(we firstly changed the order of integration and then
put $\tau = \frac{\pi}{2} - \frac{t}{2}$).
So, equating the rightmost values in the last two lines, we obtain
\be \lb{mu=}
    1-\mu_{2k}^2
 =  \frac{8}{\pi^2} \frac{4^k}{{2k \choose k}}
    \int_0^{\pi/2} t \cos^{2k}\!t\,dt\,.
\ee
Now, by Wallis inequality, we have
$$
    \sqrt{\frac{\pi}{2}} \sqrt{2k}
\le \frac{4^k}{{2k \choose k}}
\le \sqrt{\frac{\pi}{2}} \sqrt{2k+1}\,,
$$
while the integral admits the two-sided estimate
$$
      \frac{1}{2k+1}
  \le \int_0^{\pi/2} t \cos^{2k}(t)\,dt
  \le \frac{1}{2k}\,,
$$
because $\sin t \le t \le \frac{\sin t}{\cos t}$ on $[0,\frac{\pi}{2}]$,
and $\int_0^{\pi/2} \sin(t) \cos^m(t)\,dt = \frac{1}{m+1}$. Hence
$$
    \frac{8}{\pi^2}
     \sqrt{\frac{\pi}{2}}\frac{\sqrt{2k}}{2k+1}
\le 1-\mu_{2k}^2
\le \frac{8}{\pi^2} \sqrt{\frac{\pi}{2}}\frac{\sqrt{2k+1}}{2k}\,,
$$
and \rf[mu<] follows with
$c_1 = \frac{8}{\pi^2}\sqrt{\frac{\pi}{2}} \frac{2k}{2k+1} > \frac{2}{3}$
and
$c_2 = \frac{8}{\pi^2}\sqrt{\frac{\pi}{2}} \sqrt{\frac{2k+1}{2k}}
< \frac{5}{4}$.
\qed


\section{On the factor $\sqrt{r}$ at $\delta =\frac{\pi}{n}$ }

For $\delta = \frac{\pi}{n}$, our estimates for the Stechkin constant
(with the lower bound yet to be proved) look as follows:
$$
     c'\gamma_r^*
 \le K_{n,r}(\Frac{\pi}{n})
 \le c\, \sqrt{r}\ln r\, \gamma_r^*\,,
$$
i.e., the upper and lower bounds do not match. In \S\ref{1} we already
expressed our belief that additional factors on the right are
redundant. However, as we show in this section, appearance of the factor
$\sqrt{r}$ within our method is unavoidable.
(The factor $\ln r$ originates from the
use of the de la Vall\'ee Poussin sums, and perhaps can be removed by some
more sophisticated technique.)

From our initial steps \rf[W_h]-\rf[Ww], it is easy to see that
our upper estimates in all Stech\-kin inequalities
are valid  not only for the standard modulus of smoothness
$\omega_{2k}(f,h)$, but also for the modulus
\be \lb{w*}
    \omega_{2k}^*(f,h)
 := \Big\|\int_\R \wh\Delta^{2k}_t(f,\cdot) \phi_h(t)\,dt\Big\|\,,
\ee
which has a smaller value at every $h$.
It is clear that the Stechkin constant defined with respect to a smaller
modulus takes larger values, and now we show that, for the modulus
$\omega_{2k}^*(f,h)$, the increase at $h = \frac{\pi}{n}$
is exactly by the factor $\sqrt{2k}$.

\begin{theorem}
For $r=2k$, we have
$$
    \frac{\gamma_{r}^*}{1-\mu_{r}^2}
\le \sup_{f\in T_{n-1}^\perp}
    \frac{\|f\|}{\omega_{r}^*(f,\frac{\pi}{n})}
\le \frac{4}{\pi}\, \frac{\gamma_{r}^*}{1-\mu_{r}^2}\,,
$$
where
$$
       \frac{\gamma_{r}^*}{1-\mu_{r}^2}
\asymp \sqrt{r}\,\gamma_{r}^*
\asymp \frac{r}{2^{r}}\,.
$$
\end{theorem}

\proof
The upper bound was established in \rf[gBmu]-\rf[c_mu].
For the lower bound, take $f_0(x) = \cos nx$.
Then
\be \lb{cos}
    \wh\Delta^{2k}_t (f_0,x)
 = 4^k \sin^{2k}\Big(\frac{nt}{2}\Big)\cos nx, \qquad
    \phi_{\pi/n}(t)
 =  \frac{n}{\pi}\Big(1-\frac{n}{\pi}\,|t|\Big),
    \qquad |t| \le \frac{\pi}{n}\,,
\ee
hence
\baa
    \omega_{2k}^* (f_0,\Frac{\pi}{n})
&=& \Big\|\int_{-\pi/n}^{\pi/n}
     \Delta^{2k}_t(f_0,\cdot)\, \phi_{\pi/n}(t)\,dt \Big\| \\
&=& 2\cdot 4^k \int_0^{\pi/n} \sin^{2k}\Big(\frac{nt}{2}\Big)
    \frac{n}{\pi}\Big(1-\frac{n}{\pi}\,t\Big) dt \\
&=&  \frac{8}{\pi^2}\, 4^k
    \int_0^{\pi/2} \tau \cos^{2k}(\tau)\, d\tau
    \qquad \Big(\,\tau = \frac{\pi}{2} - \frac{nt}{2}\,\Big)\\
&\stackrel{\rf[mu=]}{=}& \frac{1-\mu_{2k}^2}{\gamma_{2k}^*}\,,
\eaa
while $\|f_0\| = 1$.
\qed

\medskip
Since also $E_{n-1}(f_0) = 1$, we have the same estimate
for the ratio $E_{n-1}(f_0)/\omega_{2k}^*(f_0,\frac{\pi}{n})$,
therefore, for the Stechkin constant $K_{n,r}^*(\delta)$
defined with respect to the modulus $\omega_{2k}^*(f,\delta)$,
we obtain at $\delta = \frac{\pi}{n}$
$$
    c' \sqrt{r}\, \gamma_{r}^*
\le K_{n,r}^*(\Frac{\pi}{n}) := \sup_{f\in C}
    \frac{E_{n-1}(f)}{\omega_{r}^*(f,\frac{\pi}{n})}
\le c\,\sqrt{r} \ln{r}\, \gamma_{r}^*\,.
$$


\section{Lower estimate}

\begin{lemma}
For any $n$, $r$ and $\e$, and for any $\delta < \frac{\pi}{r}$,
there exists an $f\in C$ such that,
$$
     E_{n-1}(f)
\ge \frac{1}{2} \gamma_{r-1}^*\, \omega_r(f,\delta) -\e\,.
$$
\end{lemma}

\proof
Take the step periodic function
$$
    f_0(x)
 = \left\{ \begin{array}{ll}
   1, & x \in (-\pi,0]; \\
   0, & x \in (0,\pi]. \\
   \end{array} \right.
$$
For any $x \in [-\pi,\pi]$, and for any $h < \frac{\pi}{r}$,
consider the values of this function at the points
$x_i = x+ih$, where $0 \le i \le r$.
It is clear that, for some $m \le r$, we have either
$$
    f_0(x_i) = 1, \quad 0 \le i \le m, \qquad
    f_0(x_i) = 0, \quad m < i \le r,
$$
or the other way round. Therefore, for the modulus of smoothness
$\omega_r(f_0,\delta)$, we have the following relations:
\baa
     \omega_r(f_0,\delta)
 &=& \max_{0 < h\le\delta} \max_x |\Delta_h^r f_0(x)|
  =  \max_{0 < h\le\delta} \max_x
     \Big|\sum_{i=0}^r (-1)^i {r \choose i} f_0(x+ih)\Big| \\
 &=& \max_{0 \le m \le r} \Big|\sum_{i=0}^m (-1)^i {r \choose i}\Big|
  =  \max_{0 \le m \le r} \Big|(-1)^m {r-1 \choose m}\Big|
  =  {r-1 \choose \lfloor \frac{r-1}{2}\rfloor} = 1/\gamma_{r-1}^*\,,
\eaa
i.e.,
$$
     \omega_r(f_0,\delta) = 1/\gamma_{r-1}^*\,.
$$
It is also clear that, for the best $L_\infty$-approximation of $f_0$,
we have
$$
    E_{n-1}(f_0) = \frac{1}{2},
$$
therefore the result for such an $f_0$ (without $\e$ subtracted).

This is almost what we need except that $f_0$ is not continuous.
But we can get a continuous $f$ by smoothing $f_0$ at the points
of discontinuity, say, by linearization.
For a given $\e$, set
$$
    f(x) = \frac{1}{\e} \int_{-\e/2}^{\e/2} f_0(x+t)\,dt.
$$
i.e.,
$$
    f(x)
 = \left\{ \begin{array}{l}
   1, \quad x \in [-\pi+\e,-\e]; \\
   0, \quad x \in [\e,\pi-\e]; \\
   \mbox{is linear on $[-\e,\e]$ and $[\pi-\e,\pi+\e]$.}
   \end{array} \right.
$$
Then, from the definition (or, more generally, because $f$ is
the convolution of $f_0$ with a positive kernel), it folows that
$$
    \omega_r(f,\delta) \le  \omega_r(f_0,\delta) = 1/\gamma_{r-1}^*\,.
$$
As for the best approximation of $f$, we have
$$
    E_{n-1}(f) \ge \frac{1}{2} - \e'\,.
$$
Indeed, since $E_{n-1}(f) = \|f-t_{n-1}\| \le \|f\| = 1$, the polynomial
$t_{n-1}$ of best approximation satisfies $\|t_{n-1}\| \le 2$,
therefore, by Bernstein inequality, we have
$\|t_{n-1}'\| \le 2(n-1)$,
hence, on the interval $[-\e,\e]$ of the length $2\e$
the range of $t_{n-1}$ is not more than $4(n-1)\e =: 2\e'$,
while the function $f$
on the same interval takes the values $0$ and $1$.
\qed

\begin{theorem}
For any $r$, and any $\delta \le \frac{\pi}{r}$, we have
$$
   K_{n,r}(\delta)
:= \sup_{f\in C} \frac{E_{n-1}(f)}{\omega_r(f,\delta)}
 \ge  c'_r\, \gamma_r^*
$$
where
$$
    c'_r
  = \left\{ \begin{array}{cl}
    \frac{r}{r+1}, & r = 2k-1; \\
    1,                 & r = 2k.
    \end{array} \right.
$$
In particular, for any $r$ and any $n \ge 2r$
(i.e., when $\frac{2\pi}{n} \le \frac{\pi}{r}$),
$$
   K_{n,r}\big(\Frac{2\pi}{n}\big)
:= \sup_{f\in C} \frac{E_{n-1}(f)}{\omega_r(f,\frac{2\pi}{n})}
 \ge  c'_r\, \gamma_r^*, \qquad n \ge 2r.
$$
\end{theorem}

\proof
The first lower bound is just a reformulation of the previous lemma,
because, for $\gamma_r^* := {r \choose \lfloor\frac{r}{2}\rfloor}^{-1}$,
we have $\frac{1}{2}\, \gamma_{r-1}^* = c'_r\, \gamma_r^*$\,.
\qed

\begin{remark} \rm
The order $r^{1/2} 2^{-r}$ of the lower bound for the Stechkin
constant was established earlier by Ivanov \cite{i94},
but he did not pay attention to the constant (and his extremal function
was different from ours).
\end{remark}


\section{Stechkin constants for small $r$}

For small $r=2k$, when $\mu_r$ is noticeably smaller than $1$,
our method in \S\ref{wv} will give for the Stechkin constant the
upper estimates which are better than $5\gamma_r^*$,
but they will never be smaller than $2\gamma_r^*$
because of the factor $1 + \|v_{m,n}\|$.

Surprisingly, better values (for small $r$)
which stand quite close to the lower bound
$1\cdot\gamma_r^*$ could be obtained through
technique of intermediate approximation with Steklov-type functions.
(For general $r$, this technique provides the same overblown estimate
$c_r < r^{ar}$ as Stechkin's original proof, therefore a surprise.)

Such a technique is of course well-known (it was introduced probably by
Brudnyi \cite{b64} and Freud--Popov \cite{fp70}),
and it was exploited repeatedly for proving Stechkin inequalities of various
types (e.g., for spline and one-sided approximations).
Our only innovation (if any) is the
use of the central differences instead of the forward ones,
which reduces the constants by the factor ${2k \choose k}$,
and the will to take a closer look at their actual values.

\begin{lemma}
We have
$$
    E_{n-1}(f)
\le c_{2k}\Big(\frac{\alpha\pi}{n}\Big)\,
    \gamma_{2k}^*\, \omega_{2k}\Big(f,\frac{\alpha\pi}{n}\Big)\,,
$$
where
\be \lb{ca1}
   c_{2k}\Big(\frac{\alpha\pi}{n}\Big)
 = 1 + F_{2k}\frac{k^{2k}}{(\alpha\pi)^{2k}}
   \sum_{i=1}^k \frac{2b_i}{i^{2k}}, \qquad b_i = {2k \choose k+i}\,,
\ee
and $F_{2k}$ are the Favard constants.
\end{lemma}

\proof
Given $f$, with any $2k$ times differentiable function $f_h$,
we have
\be \lb{sm}
      E_{n-1}(f)
 \;\le\; E_{n-1}(f-f_h) + E_{n-1}(f_h)
 \;\le\;  \|f-f_h\| + \frac{F_{2k}}{n^{2k}}\,\|f_h^{(2k)}\|
\ee
where we used the Favard inequality for the best approximations of $f_h$.
A typical choice of $f_h$ is via the Steklov functions of order $2k$:
\baa
& \displaystyle
     I_{ih}(f,x)
 :=  \frac{1}{(h/k)^{2k}}
     \underbrace{\int_{-h/2k}^{h/2k}\cdots\int_{-h/2k}^{h/2k}}_{2k}
     f(x- i(t_1\!+\!\cdots\!+\!t_{2k}))\,dt_1\cdots dt_{2k}\,, & \\
& \displaystyle
   I_{ih}^{(2k)}(f,x)
 = \frac{(-1)^k}{(ih/k)^{2k}}\wh\Delta^{2k}_{ih/k}f(x)\,, &
\eaa
namely
$$
     f_h
  := \frac{1}{{2k\choose k}} \sum_{{i=-k \atop i\ne 0}}^{k} (-1)^{i+1}
     {2k\choose k+i} I_{ih}(f)
   = \gamma_{2k}^* \sum_{i=1}^{k} (-1)^{i+1} 2b_i I_{ih}(f)\,.
$$
Then
\baa
&       \|f - f_h\|
 \le \gamma_{2k}^*\,\omega_{2k}(f,h), & \\
&  \displaystyle     \|f_h^{(2k)}\|
 \le \gamma_{2k}^*\,\sum_{i=1}^k \frac{2b_i}{(ih/k)^{2k}}\,
       \omega_{2k}(f,ih/k)
 \;\le\;  \gamma_{2k}^*\,\omega_{2k}(f,h)\,
      \frac{k^{2k}}{h^{2k}} \sum_{i=1}^k \frac{2b_i}{i^{2k}}\,,&
\eaa
whence applying \rf[sm]
$$
      E_{n-1}(f)
\;\le\; c_{2k}(h)\, \gamma_{2k}^*\,\omega_2(f,h)\,,\qquad
   c_{2k}(h)
 = 1 + F_{2k}\frac{k^{2k}}{(nh)^{2k}}
   \sum_{i=1}^k \frac{2b_i}{i^{2k}}\,,
$$
and we take $h = \frac{\alpha\pi}{n}$.
\qed

\medskip
In \rf[ca1], we can obtain a small value only if
$\frac{k}{\alpha\pi} < 1$,
i.e., we may try $k = (1,2,3)$ for $\alpha = 1$, and
$k = (1,2,3,4,5)$ for $\alpha = 2$. So we did (dropping
those values for which the resulting constants in \rf[ca1] were not
close to $1$).

\begin{theorem} \lb{small r}
For $\delta = \frac{\pi}{n}$ and $\delta = \frac{2\pi}{n}$, we have
$$
    E_{n-1}(f)
\le c_{r}(\delta)\,\gamma_{r}^*\, \omega_{r}(f,\delta)\,,
$$
where $c_{2k-1}(\delta) = c_{2k}(\delta)$, and the values of
$c_{2k}(\delta)$ are given below
$$
\begin{array}{c|c}
c_2(\frac{\pi}{n}) & c_4(\frac{\pi}{n}) \\  \hline
1\frac{1}{4}       & 2\frac{7}{10}
\end{array}\,,
\qquad
\begin{array}{c|c|c}
c_2(\frac{2\pi}{n}) & c_4(\frac{2\pi}{n}) & c_6(\frac{2\pi}{n}) \\  \hline
     1\frac{1}{16}  & 1\frac{1}{9}        & 1\frac{1}{2}
\end{array}\,.
$$
\end{theorem}

\proof
We will use the following values:
$F_2 = \frac{\pi^2}{8}$,
$F_4 = \frac{5\pi^4}{384}$,
$F_6 = \frac{61\pi^6}{46080}$.

1) For $2k=2$, we have
$$
    c_{2}\Big(\frac{\alpha\pi}{n}\Big)
  = 1 + \frac{\pi^2}{8} \frac{2}{(\alpha\pi)^2}
  = 1 + \frac{1}{4\alpha^2}\,.
$$
With $\alpha=1$ and $\alpha=2$, we obtain
$c_{2}(\frac{\pi}{n}) = \Frac{5}{4}$ and
$c_{2}(\frac{2\pi}{n}) = \Frac{17}{16}$.
Also, with $\alpha=\frac{1}{2}$, we obtain the remarkable inequality
$$
   E_{n-1}(f) \le 1\cdot \omega_2\Big(f,\frac{\pi}{2n}\Big)\,.
$$

2) For $2k=4$,
$$
    c_{4}\Big(\frac{\alpha\pi}{n}\Big)
  = 1 + \frac{5\pi^4}{384} \frac{2^4}{(\alpha\pi)^4} \cdot
    2 \left[\frac{4}{1^4} + \frac{1}{2^4}\right]
  = 1 + \frac{325}{192} \frac{1}{\alpha^4}\,.
$$
With $\alpha=1$ and $\alpha=2$, we obtain
$c_{4}(\frac{\pi}{n}) = \Frac{517}{192} = 2.6927$, and
$c_{4}(\frac{2\pi}{n}) = \Frac{3397}{3072} = 1.1058$.

\medskip
3) For $2k=6$, with $\alpha=2$, we have
$$
     c_{6}\Big(\frac{2\pi}{n}\Big)
  =  1 + \frac{61\pi^6}{46080}\,\frac{3^6}{(2\pi)^6}\cdot
     2 \left[\frac{15}{1^4} + \frac{6}{2^6} + \frac{1}{3^6}\right]
  = 1.4552 < 1\frac{1}{2}\,.
\eqno{\Box}
$$

Theorem \ref{small r} provides a certain support to our Conjecture \ref{con},
which says, in particular,
that, for even $r=2k$, and for $\delta \ge \frac{\pi}{n}$,
the best constant in the Stechkin inequality
has the value $K_{n,r}(\delta) = 1\cdot\gamma_r^*$.

\medskip
{\bf Acknowledgements.}
Our thanks to Alexander Babenko
for his comments on a draft of this paper.

\small

\begin{tabular}[t]{l}
Simon Foucart \\
Dept of Mathematics \\
1533 Stevenson Center \\
Vanderbilt University \\
Nashville, TN 37240 \\
USA \\
\textrm{simon.foucart@vanderbilt.edu}
\end{tabular}
\hfill
\begin{tabular}[t]{l}
Yuri Kryakin \\
Institute of Mathematics \\
University of Wroclaw \\
Plac Grunwaldzki 2/4 \\
50-384 Wroclaw \\
Poland \\
\textrm{kryakin@list.ru}
\end{tabular}
\hfill
\begin{tabular}[t]{l}
Alexei Shadrin \\
DAMTP \\
University of Cambridge \\
Wilberforce Road \\
Cambridge CB3 0WA \\
UK \\
\textrm{a.shadrin@damtp.cam.ac.uk}
\end{tabular}

\end{document}